\documentclass[12pt,a4paper,reqno]{amsart}
\usepackage{amsmath}
\usepackage{amsfonts}
\usepackage{amssymb}
\usepackage{bm}
\usepackage{graphics}

\renewcommand\Re{{\operatorname{Re}}}
\renewcommand\Im{{\operatorname{Im}}}
\newcommand\Vol{\hbox{\rm Vol}}

\newcommand\R{{\mathbf{R}}}
\newcommand\Z{{\mathbf{Z}}}
\newcommand\BZ{{\mathbf{Z}}}
\newcommand\BBC{{\mathbf{C}}}
\newcommand\C{{\mathbf{C}}}
\newcommand\Q{{\mathbf{Q}}}

\renewcommand\P{{\mathbf{P}}}
\newcommand\E{{\mathbf{E}}}

\newcommand\eps{\varepsilon}

\newcommand\trace{{\operatorname{trace}}}
\newcommand\tr{{\operatorname{trace}}}

\renewcommand\a{{\operatorname{x}}}
\renewcommand\b{{\operatorname{y}}}
\newcommand\bv{{\mathbf{v}}}

\theoremstyle{plain}
  \newtheorem{theorem}[subsection]{Theorem}
  \newtheorem{conjecture}[subsection]{Conjecture}

  \newtheorem{lemma}[subsection]{Lemma}
  \newtheorem{cor}[subsection]{Corollary}

\theoremstyle{remark}
  \newtheorem{remark}[subsection]{Remark}
  
  \newtheorem{example}[subsection]{Example}

\theoremstyle{definition}

\include{psfig}

\begin{document}

\title[Littlewood-Offord, circular law, universality]{ From the Littlewood-Offord problem to the Circular Law: universality of the spectral distribution of random matrices }

\author{Terence Tao}
\address{Department of Mathematics, UCLA, Los Angeles CA 90095-1555}
\email{tao@math.ucla.edu}
\thanks{T. Tao is supported by NSF grant CCF-0649473 and a grant from the MacArthur Foundation.}

\author{Van Vu}
\address{Department of Mathematics, Rutgers, Piscataway, NJ 08854}
\email{vanvu@math.rutgers.edu}

\thanks{V. Vu is  is supported by  NSF Career Grant 0635606.}


\begin{abstract}  The famous \emph{circular law}  
asserts that if $M_n$ is an $n \times n$ matrix with iid complex entries of mean zero and unit variance, then the empirical spectral distribution (ESD) of the normalized matrix $\frac{1}{\sqrt{n}} M_n$ converges both in probability and almost surely to the uniform distribution on the unit disk $\{ z \in \C: |z| \leq 1 \}$.  After a long sequence of partial results that verified this law under additional assumptions on the distribution of the entries, the circular law is now known to be true for arbitrary distributions with mean zero and unit variance.  In this survey we describe some of the key ingredients used in the establishment of the circular law at this level of generality, in particular recent advances in understanding the Littlewood-Offord problem and its inverse.
\end{abstract}

\maketitle

\section {ESD of random matrices}

For an $n \times n$ matrix $A_n$ with complex entries, let
$$\mu_{A_n} := \frac{1}{n} \sum_{i=1}^n \delta_{\lambda_i}$$
be the \emph{empirical spectral distribution} (ESD) of its
eigenvalues $\lambda_i \in \BBC, i=1, \dots n$ (counting multiplicity), thus for instance
$$ \mu_{A_n}( \{ z \in \BBC | \Re z \leq s; \Im z \leq t \} ) = \frac{1}{n} | \{ 1 \leq i \leq n: \Re \lambda_i \leq s; \Im \lambda_i \leq t \}|$$
for any $s,t \in \R$ (we use $|A|$ to denote the cardinality of a finite set $A$), and
$$ \int f\ d\mu_{A_n} = \frac{1}{n} \sum_{i=1}^n f(\lambda_i)$$
for any continuous compactly supported $f$. Clearly, $\mu_{A_n}$ is a discrete probability measure on $\BBC$.

A fundamental problem in the theory of random matrices is to compute
the limiting distribution of the ESD $\mu_{A_n}$ of a sequence of random
matrices $A_n$ with sizes tending to infinity \cite{Mehta, BS}. In what
follows, we consider normalized random matrices of the form $A_n = \frac{1}{\sqrt{n}} M_n$, where $M_n = (\a_{ij})_{1 \leq i,j \leq n}$ has entries that are iid random variables $\a_{ij} \equiv \a$.  Such matrices have been studied at least as far back as Wishart \cite{wish} (see \cite{Mehta, BS} for more discussion).

One of the first limiting distribution results is the famous semi-circle law
of Wigner \cite{wig}. Motivated by research in nuclear physics,
Wigner studied Hermitian random matrices with (upper triangular)
entries being iid random variables with mean zero and variance one.  In the Hermitian case, of course, the ESD is supported on the real line $\R$.  He proved that the expected ESD of a normalized $n \times n$ Hermitian matrix $\frac{1}{\sqrt{n}} M_n$, where $M_n = (\a_{ij})_{1 \leq i,j \leq n}$ has iid gaussian entries $\a_{ij} \equiv N(0,1)$, converges in the sense of probability measures\footnote{We say that a collection $\mu_n$ of probability measures converges to a limit $\mu$ if one has $\int f\ d\mu_n \to \int f\ d\mu$ for every continuous compactly supported function $f$, or equivalently if $\mu( \{ z \in \BBC | \Re z \leq s; \Im z \leq t \} )$ converges to $\mu( \{ z \in \BBC | \Re z \leq s; \Im z \leq t \} )$ for all $s, t$.} to the
semi-circle distribution
\begin{equation}\label{semicircle}
 \frac{1}{2\pi} 1_{[-2,2]}(x) \sqrt {4-x^2}\ dx
\end{equation}
on the real line, where $1_E$ denotes the indicator function of a set $E$.

\begin{theorem}[Semi-circular law for the Gaussian ensemble]\label{theorem:Wigner}\cite{wig}
Let $M_n$ be an $n \times n$ random Hermitian matrix whose  entries
are iid gaussian variables with mean 0 and variance 1. Then, with
probability one, the ESD of $\frac{1}{\sqrt n} M_n$ converges in the sense of probability measures to the
semi-circle law \eqref{semicircle}.
\end{theorem}

Henceforth we shall say that a sequence $\mu_n$ of random probability measures converges \emph{strongly} to a deterministic probability measure $\mu$ if, with probability one, $\mu_n$ converges in the sense of probability measures to $\mu$.  We also say that $\mu_n$ converges \emph{weakly} to $\mu$ if for every continuous compactly supported $f$, $\int f\ d\mu_n$ converges in probability to $\int f\ d\mu$, thus $\P( |\int f\ d\mu_n - \int f\ d\mu| > \eps ) \to 0$ as $n \to \infty$ for each $\eps > 0$.  Of course, strong convergence implies weak convergence; thus for instance in Theorem \ref{theorem:Wigner}, $\mu_{\frac{1}{\sqrt{n}} M_n}$ converges both weakly and strongly to the semicircle law.

Wigner also proved similar results for various other distributions,
such as the Bernoulli distribution (in which each $\a_{ij}$ equals $+1$ with probability $1/2$ and $-1$ with probability $1/2$). His work has been extended and strengthened in several aspects \cite{Arnold1, Arnold2, Pastur}. The most general form was proved by Pastur \cite{Pastur}:

\begin{theorem}[Semi-circular law]\label{theorem:Pastur} \cite{Pastur}
Let $M_n$ be an $n \times n$ random Hermitian matrix whose entries
are iid complex random variables with mean 0 and variance 1. Then ESD of $\frac{1}{\sqrt n} M_n$ converges (in both the strong and weak senses) to
the semi-circle law.
\end{theorem}

The situation with non-Hermitian matrices is much more complicated, due to the presence of \emph{pseudospectrum}\footnote{Informally, we say that a complex number $z$ lies in the pseudospectrum of a square matrix $A$ if $(A-zI)^{-1}$ is large (or undefined).  If $z$ lies in the pseudospectrum, then small perturbations of $A$ can potentially cause $z$ to fall into the spectrum of $A$, even if it is initially far away from this spectrum.  Thus, whenever one has pseudospectrum far away from the actual spectrum, the actual distribution of eigenvalues can depend very sensitively (in the worst case) on the coefficients of $A$.  Of course, our matrices are random rather than worst-case, and so we expect the most dangerous effects of pseudospectrum to be avoided; but this of course requires some analytical effort to establish, and deterministic techniques (e.g. truncation) should be used with extreme caution, since they are likely to break down in the worst case.}  that can potentially make the ESD quite unstable with respect to perturbations.  The non-Hermitian variant of this theorem, the Circular Law
Conjecture, has been raised since the 1950's  (see Chapter 10 of \cite{BS} or the introduction of \cite{bai})

\begin{conjecture}[Circular law]\label{conj:CL}
Let $M_n$ be the $n \times n$  random matrix whose  entries are iid
complex random variables with mean 0 and variance 1. Then the ESD of
$\frac{1}{\sqrt n} M_n$ converges (in both the strong and weak
senses) to the uniform distribution $\mu := \frac{1}{\pi} 1_{|z| \leq 1} dz$ on the unit disk $\{ z \in \C: |z| \leq 1 \}$.
\end{conjecture}

The numerical evidence for this conjecture is extremely strong (see e.g. Figure \ref{figure:CircLaw}).   However, there are significant difficulties in establishing this conjecture rigorously, not least of which is the fact that the main techniques used to handle Hermitian matrices (such as moment methods and truncation) can not be applied to the non-Hermitian model (see \cite[Chapter 10]{BS} for a detailed discussion).  Nevertheless, the conjecture has been intensively worked on for many decades.  The circular law was verified for the complex gaussian distribution in \cite{Mehta} and the real gaussian distribution in \cite{edel}.  An approach to attack the general case was introduced in \cite{Girko1}, leading to a resolution of the strong circular law for continuous distributions with bounded sixth moment in \cite{bai}.  The sixth moment hypothesis in \cite{bai} was lowered to $(2+\eta)^{\operatorname{th}}$ moment for any $\eta > 0$ in \cite{BS}.  The removal of the hypothesis of continuous distribution required some new ideas.  In \cite{GT1} the weak circular law for (possibly discrete) distributions with subgaussian moment was established, with the subgaussian condition relaxed to a fourth moment condition in \cite{PZ} (see also \cite{Girko2} for an earlier result of similar nature), and then to $(2+\eta)^{\operatorname{th}}$ moment in \cite{GT2}.  Shortly before this last result, the strong circular law assuming $(2+\eta)^{\operatorname{th}}$ moment was established in \cite{TVcir1}.  Finally, in a recent paper \cite{TVcir2}, the authors proved this conjecture (in both strong and weak forms) in full generality. In fact, we obtained this result as a consequence of a more general theorem, presented in the next section.

\section{Universality}

An easy case of Conjecture \ref{conj:CL} is when the entries $\a_{ij}$ of $M_n$ are iid
complex gaussian. In this case there is the following precise formula
 for the joint density function of the eigenvalues, due to Ginibre \cite{gin} (see also \cite{Mehta}, \cite{hwang} for more discussion of this formula):

\begin{equation} \label{eqn:Ginibre} p(\lambda_{1}, \cdots, \lambda_{n}) = c_{n} \prod_{[i <
j}|\lambda_{i}- \lambda_{j}|^{2 } \prod_{i=1} ^{n} e^{-n
|\lambda_{i}|^{2}}. \end{equation}

From here one can verify the conjecture in this case by a direct
calculation. This was first done by Mehta and also Silverstein in the
1960s:

\begin{theorem}[Circular law for Gaussian matrices]\label{theorem:mehta} \cite{Mehta}
Let $M_n$ be an $n \times n$ random  matrix whose  entries are iid
complex gaussian variables with mean 0 and variance 1. Then, with
probability one, the ESD of $\frac{1}{\sqrt n} M_n$ tends to the
circular law. \end{theorem}

A similar result for the real gaussian ensemble was established in \cite{edel}.  These methods rely heavily on the strong symmetry properties of such ensembles (in particular, the invariance of such ensembles with respect to large matrix groups such as $O(n)$ or $U(n)$) in order to perform explicit algebraic computations, and do not extend directly to more combinatorial ensembles, such as the Bernoulli ensemble.

The above mentioned results and conjectures can be viewed as examples of
a general phenomenon in probablity and mathematical physics, namely, that
global information about a large random system (such as
limiting distributions) does not depend on the particular
distribution of the particles. This is often referred to as the
\emph{universality} phenomenon (see e.g. \cite{deift}). The most famous example of this
phenomenon is perhaps the central limit theorem.

In view of the universality phenomenon, one can
see that Conjecture \ref{conj:CL} generalizes  Theorem
\ref{theorem:mehta} in the same way that Theorem \ref{theorem:Pastur}
generalizes  Theorem \ref{theorem:Wigner}.

\begin{figure}
\centerline{\textbf{Bernoulli \hspace{1.8in} Gaussian }}
\begin{center}
\scalebox{.3}{\includegraphics{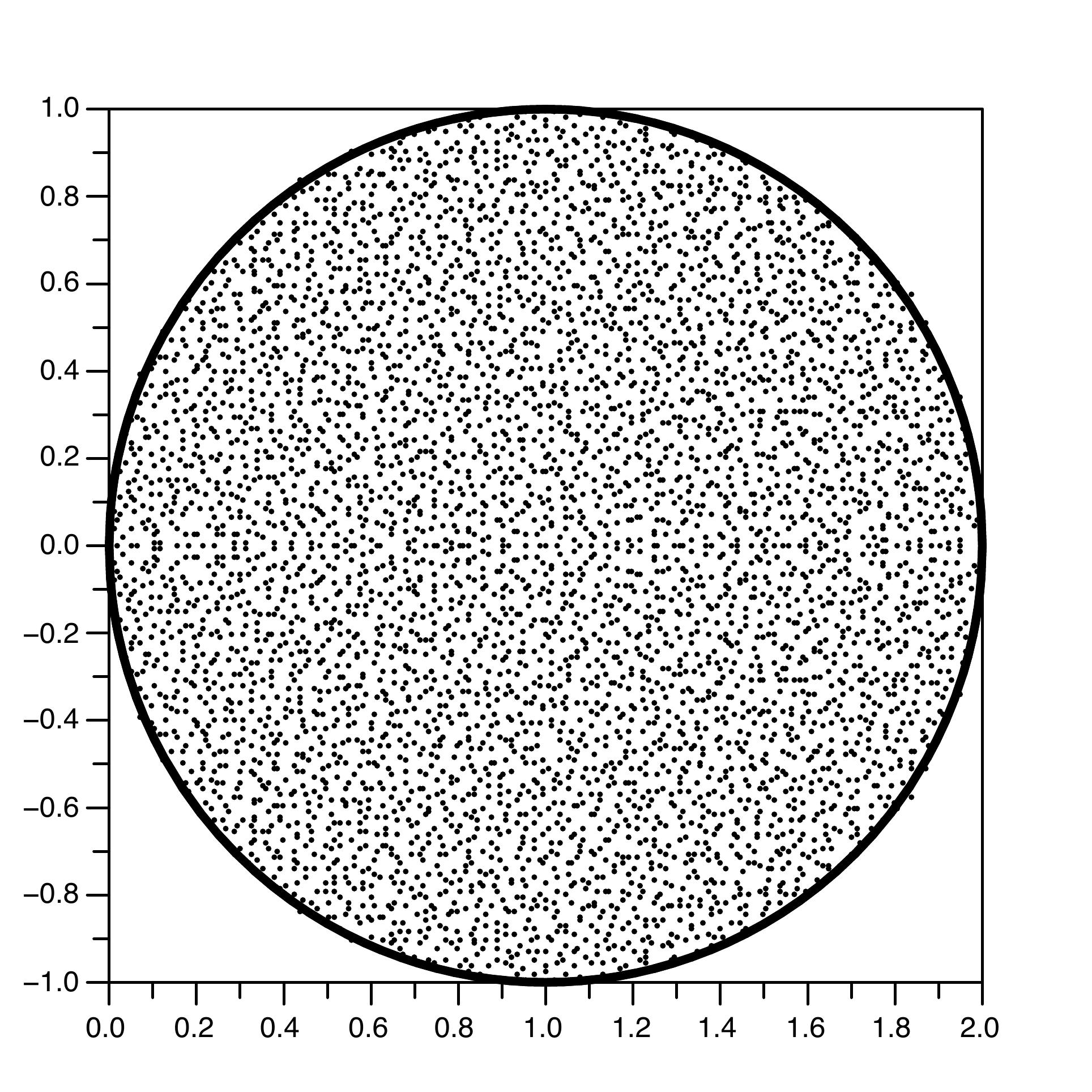}}
\scalebox{.3}{\includegraphics{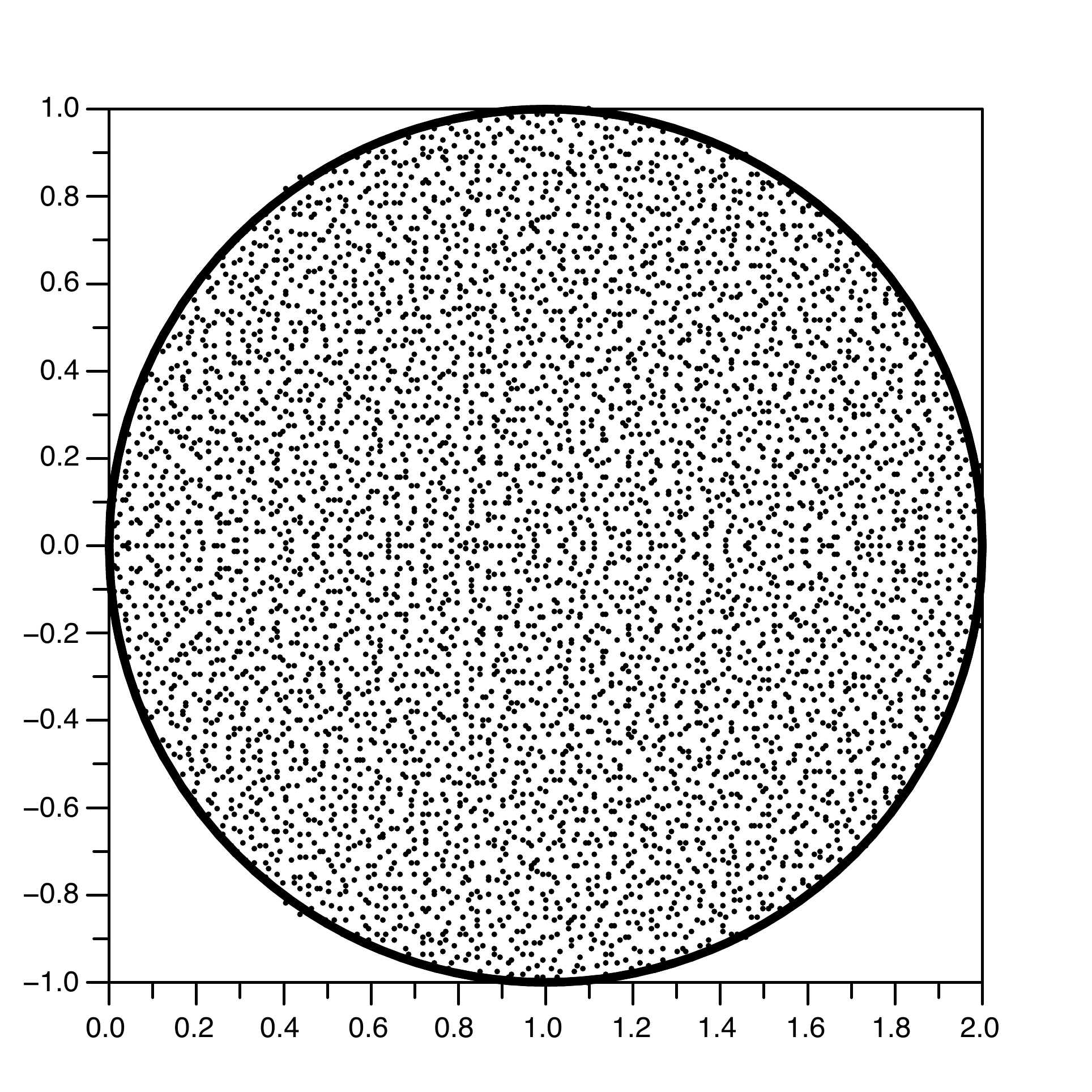}}
\end{center}
\caption{Eigenvalue plots of two randomly generated 5000 by 5000 matrices.  On the left, each entry was an iid Bernoulli random variable, taking the values $+1$ and $-1$ each with probability $1/2$.  On the right, each entry was an iid Gaussian normal random variable, with probability density function is $\frac{1}{\sqrt{2\pi}} \exp( -x^2/2 ) $.  (These two distributions were shifted by adding the identity matrix, thus the circles are centered at $(1,0)$ rather than at the origin.)}
\label{figure:CircLaw}
\end{figure}

A demonstration of the circular law for the Bernoulli and the
Gaussian case appears\footnote{We thank Phillip Wood for creating the figures in this paper.} in the Figure~\ref{figure:CircLaw}.

\begin{figure}
\centerline{\textbf{Bernoulli \hspace{1.8in} Gaussian }}
\begin{tabbing}
\qquad
\=
\scalebox{.3}{\includegraphics{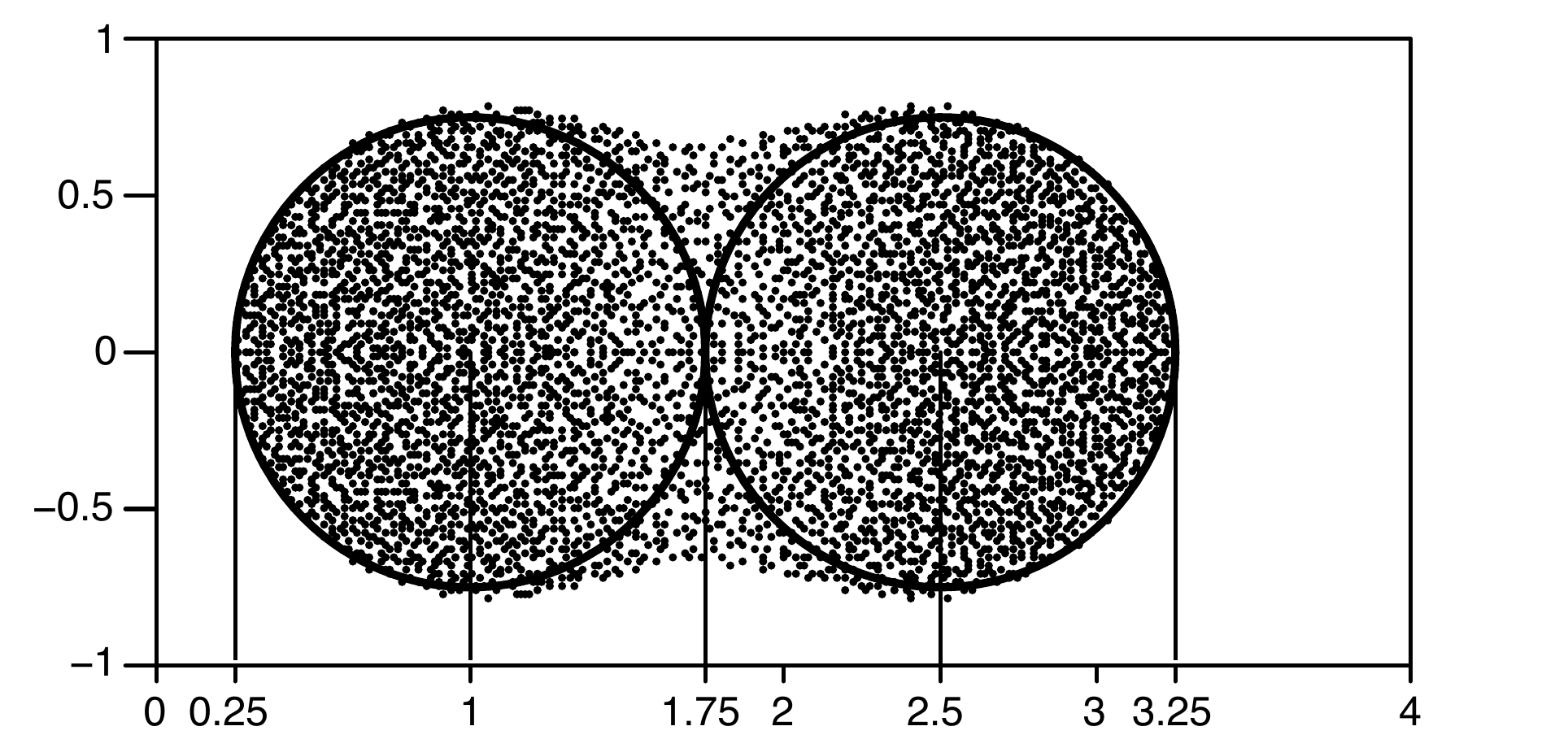}}
\=
\scalebox{.3}{\includegraphics{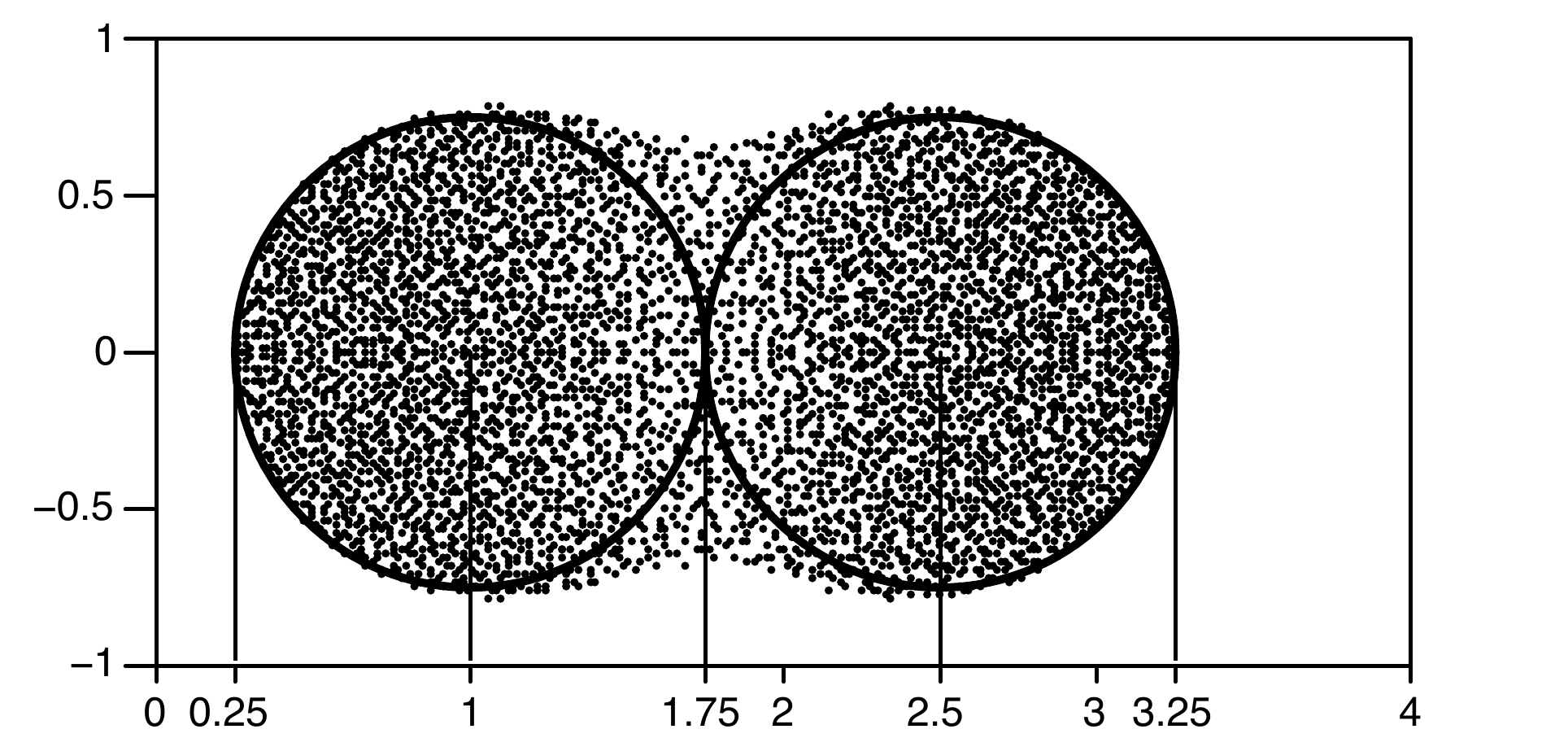}} \\
\>
\scalebox{.3}{\includegraphics{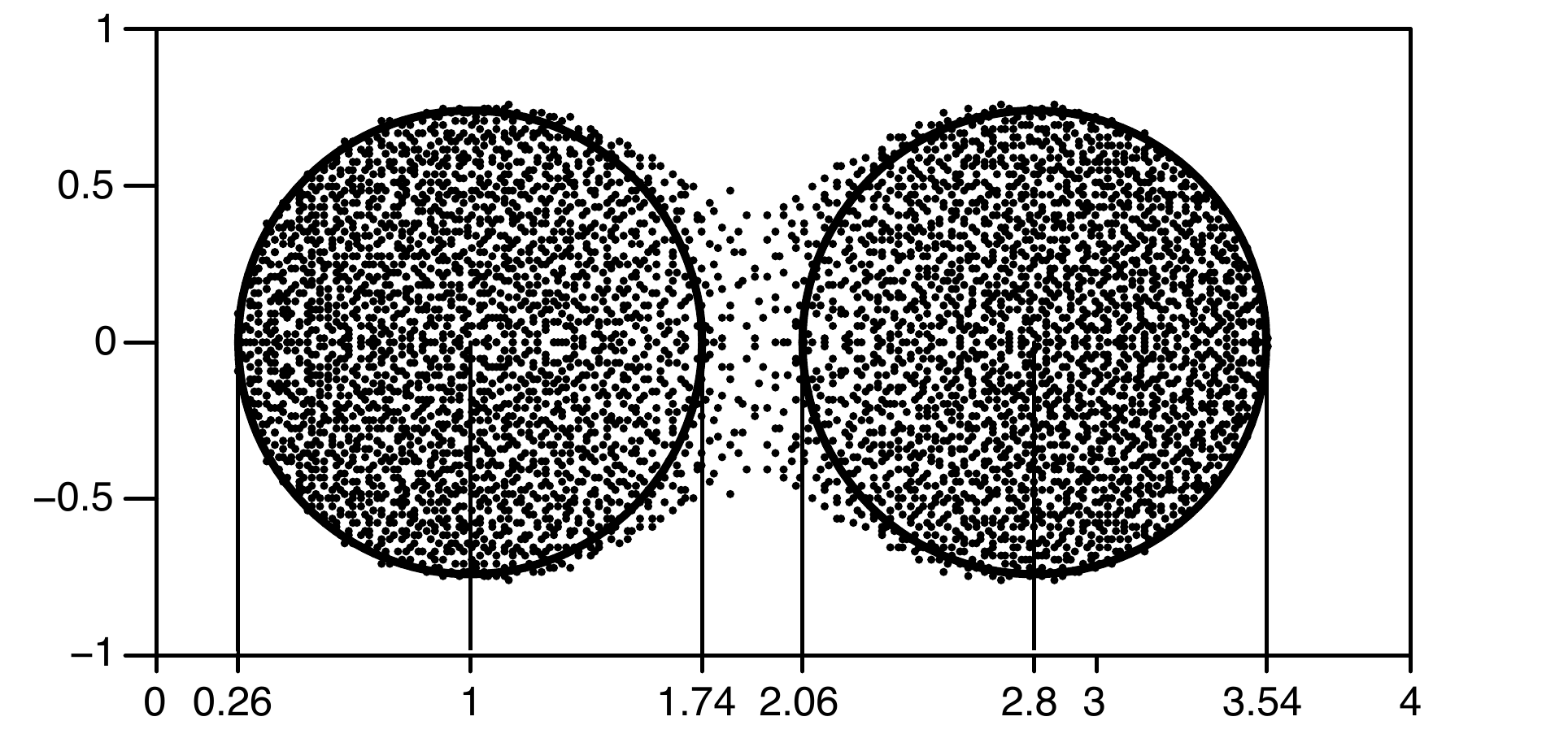}}
\>
\scalebox{.3}{\includegraphics{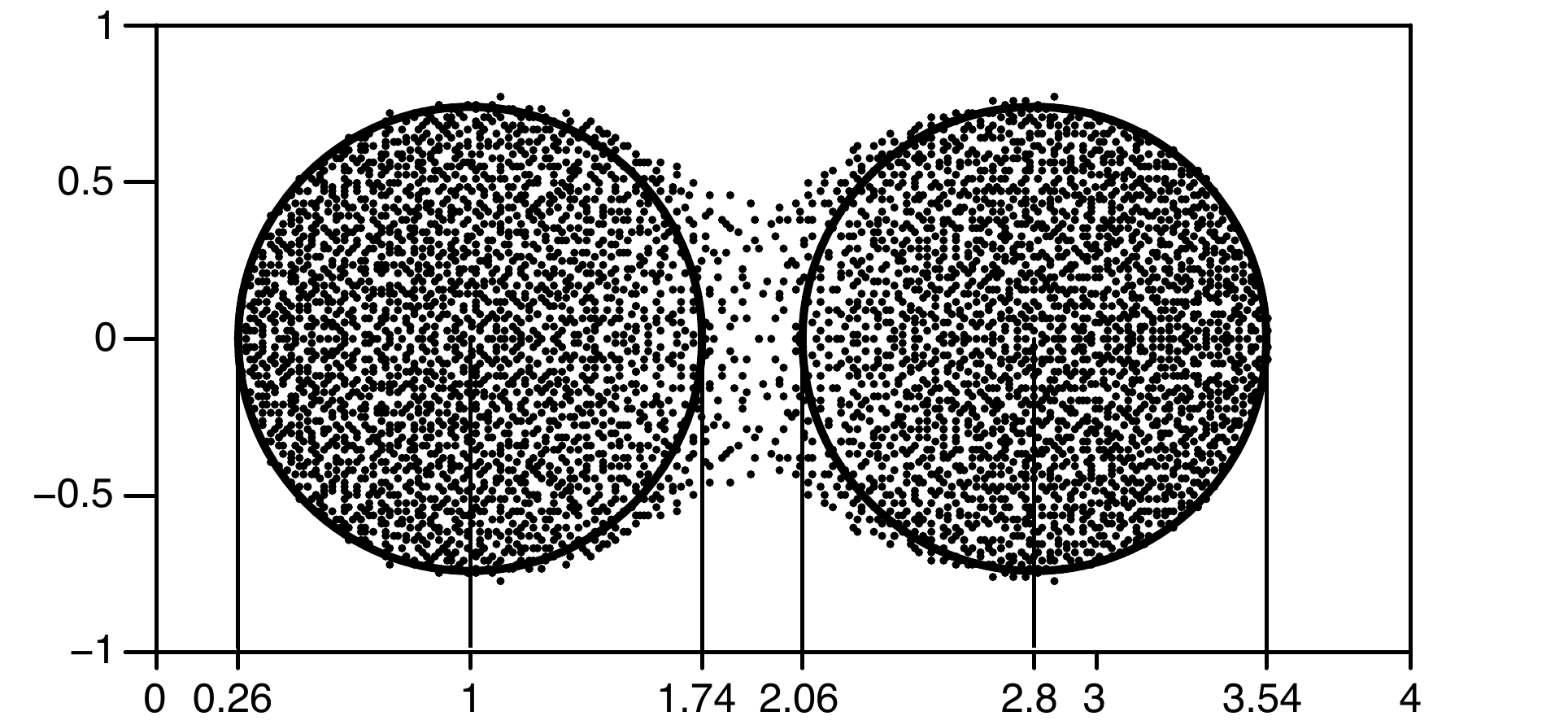}}
\end{tabbing}
\caption{Eigenvalue plots of randomly generated $n$ by $n$ matrices of the form $D_n+M_n$, where $n=5000$.  In the left column, each entry of $M_n$ was an iid Bernoulli random variable, taking the values $+1$ and $-1$ each with probability $1/2$, and in the right column, each entry was an iid Gaussian normal random variable, with probability density function is $\frac{1}{\sqrt{2\pi}} \exp( -x^2/2 )$.  In the first row, $D_n$ is the deterministic matrix $\operatorname{diag}(1,1,\ldots,1,2.5,2.5,\ldots,2.5)$, and in the second row $D_n$ is the deterministic matrix $\operatorname{diag}(1,1,\ldots,1,2.8,2.8,\ldots,2.8)$ (in each case, the first $n/2$ diagonal entries are $1$'s, and the remaining entries are $2.5$ or $2.8$ as specified). }
\label{figure:ChangeMean}
\end{figure}

\begin{figure}
\centerline{\textbf{Bernoulli \hspace{1.8in} Gaussian }}
\begin{center}
\scalebox{.32}{\includegraphics{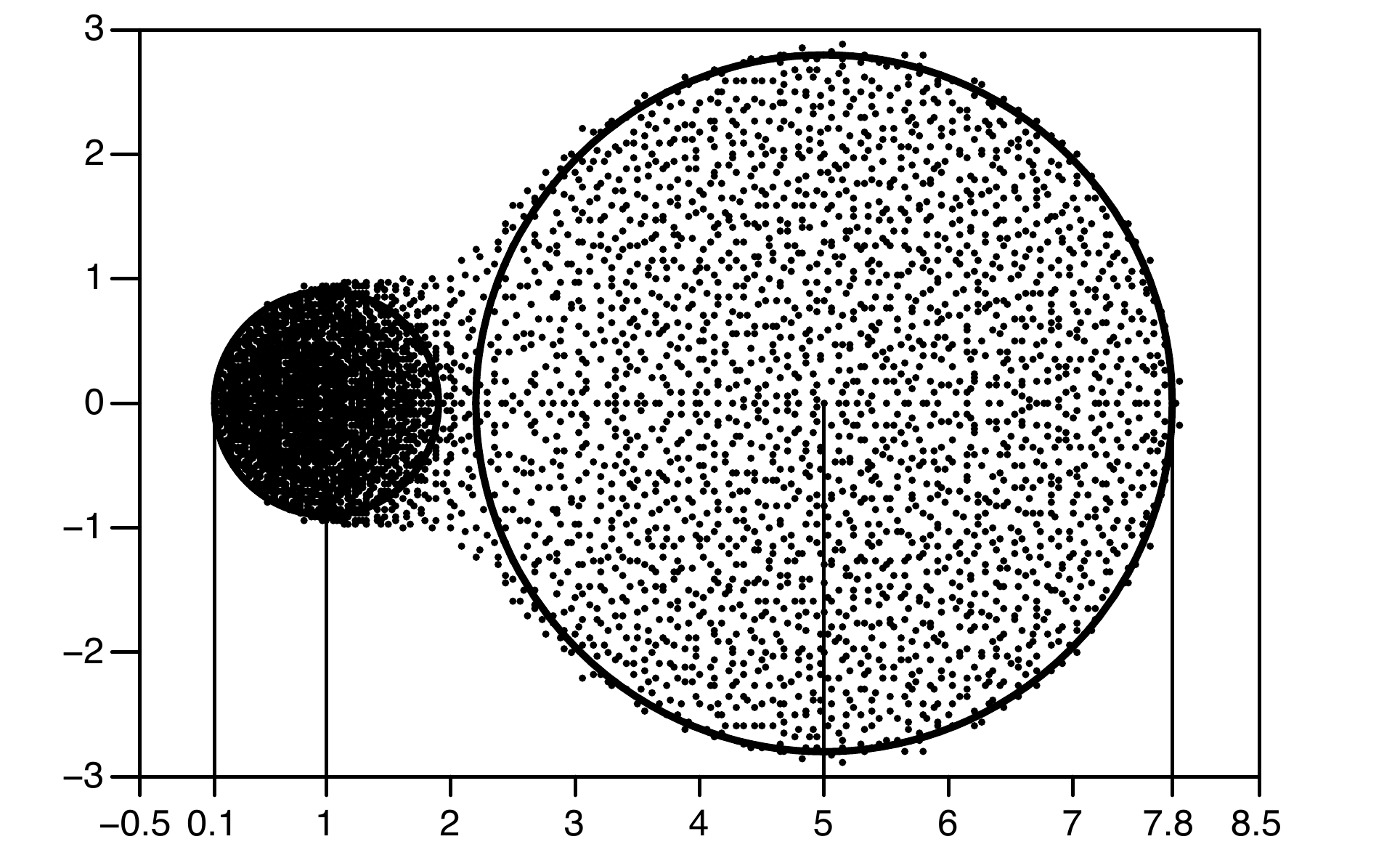}}
\scalebox{.32}{\includegraphics{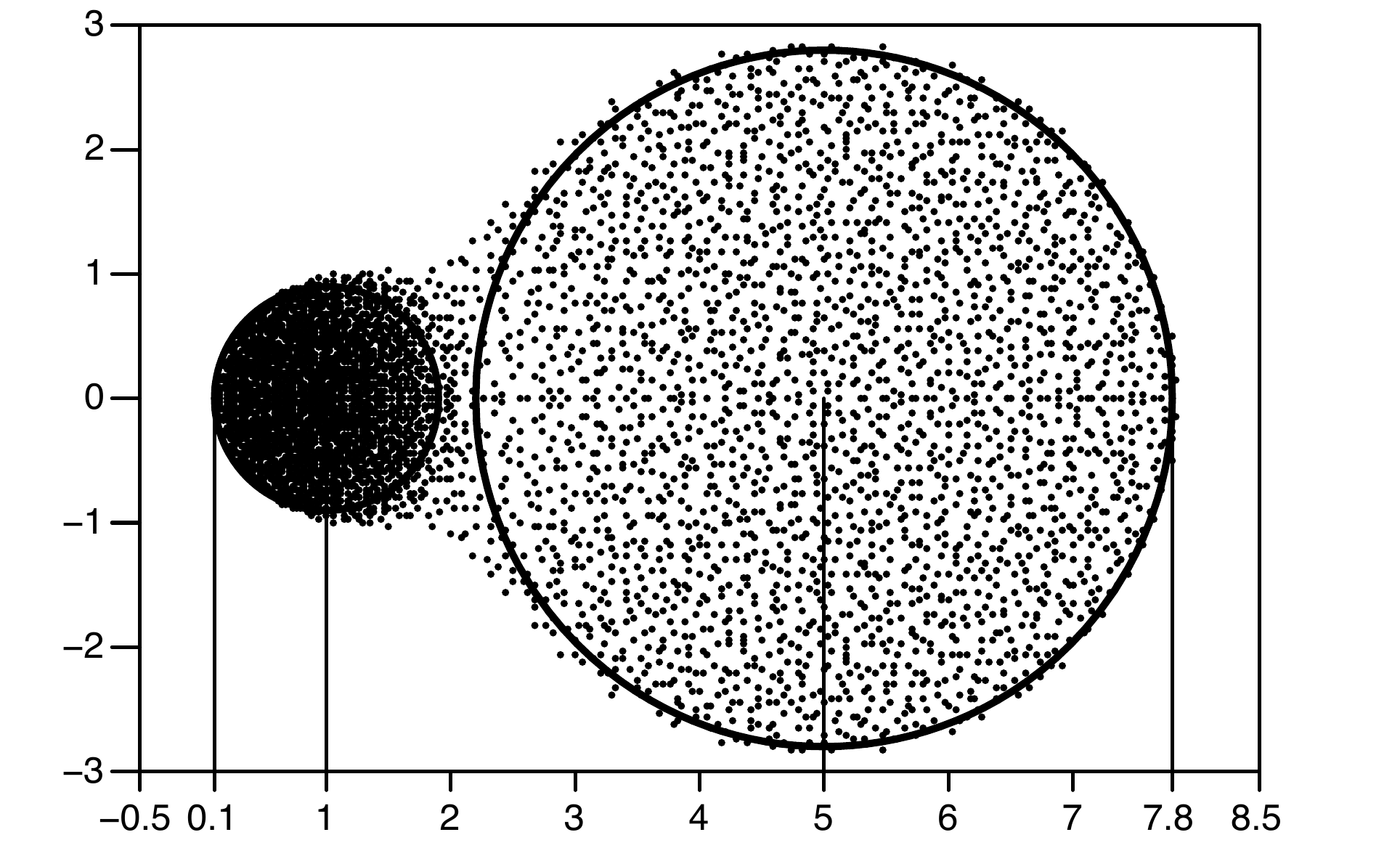}}
\end{center}
\caption{Eigenvalue plots of two randomly generated 5000 by 5000 matrices of the form $A + BM_nB$, where $A$ and $B$ are diagonal matrices having $n/2$ entries with the value 1 followed by $n/2$ entries with the value 5 (for $D$) and the value $2$ (for $X$).
On the left, each entry of $M_n$ was an iid Bernoulli random variable, taking the values $+1$ and $-1$ each with probability $1/2$.  On the right, each entry of $M_n$ was an iid Gaussian normal random variable, with probability density function is $\frac{1}{\sqrt{2\pi}} \exp( -x^2/2 ) $.}
\label{figure:Extension}
\end{figure}

The universality phenomenon seems to hold even for more general
models of random matrices, as demonstrated by Figure~\ref{figure:ChangeMean} and Figure~\ref{figure:Extension}.

This evidence  suggests that the asymptotic shape of the ESD depends
only on the mean and the variance of each entry in the matirx. As mentioend earlier, the main result
of \cite{TVcir2} (building on a large number of previous results) gives a rigorous proof of this phenomenon in full generality.

For any matrix $A$, we define the \emph{Frobenius norm} (or \emph{Hilbert-Schmidt norm})
$\|A\|_F$ by the formula $\|A\|_F := \trace(AA^\ast)^{1/2} =
\trace(A^\ast A)^{1/2}$.

\begin{theorem}[Universality principle]\label{theorem:main1}
Let $\a$ and $\b$ be complex random variables with zero mean and
unit variance. Let $X_n =  (\a_{ij})_{1 \leq i,j \leq n}$ and $Y_n
:= (\b_{ij})_{1 \leq i,j \leq n}$ be $n \times n$ random matrices
whose entries $\a_{ij}$, $\b_{ij}$ are iid copies of $\a$ and $\b$,
respectively. For each $n$, let $M_n$ be a deterministic $n \times
n$ matrix satisfying
\begin{equation} \label{eqn:conditionM}  \sup_n \frac{1}{n^2} \|M_n\|_F^2 < \infty. \end{equation}
Let $A_n:= M_n + X_n$ and $B_n:= M_n +Y_n$. Then
$\mu_{\frac{1}{\sqrt{n}} A_n} - \mu_{\frac{1}{\sqrt{n}} B_n}$
converges weakly to zero.  If furthermore we make the additional
hypothesis that the ESDs \begin{equation}\label{must}
\mu_{(\frac{1}{\sqrt{n}} M_n-zI) (\frac{1}{\sqrt{n}} M_n-zI)^\ast}
\end{equation}
converge in the sense of probability measures to a limit for almost every $z$, then
$\mu_{\frac{1}{\sqrt{n}} A_n} - \mu_{\frac{1}{\sqrt{n}} B_n}$
converges strongly to zero.
\end{theorem}

This theorem reduces the computing of the limiting distribution to
the case where one can assume\footnote{Some related ideas also appear in \cite{Girko2}.  In the context of the central limit theorem, the idea of replacing arbitrary iid ensembles by Gaussian ones goes back to Lindeberg \cite{Lind}, and is sometimes known as the \emph{Lindeberg invariance principle}; see \cite{chat} for further discussion, and a formulation of this principle for Hermitian random matrices.}
 that the entries $\a$ have Gaussian (or any
special) distribution. Combining this theorem (in the case $M_n = 0$) with
Theorem \ref{theorem:mehta}, we conclude

\begin{cor}
The circular law (Conjecture \ref{conj:CL}) holds in both the weak and strong sense.
\end{cor}

It is useful to notice that Theorem \ref{theorem:main1} still holds
even when the limiting distributions do not exist.

The proof of Theorem \ref{theorem:main1} relies on several
surprising connections between seemingly remote areas of
mathematics that have been discovered in the last few years. The
goal of this article is to give the reader an overview of these
connections and through them a sketch of the proof of Theorem
\ref{theorem:main1}. The first area we shall visit is combinatorics.

\section {Combinatorics}

As we shall discuss later, one of the primary difficulties in controlling the ESD of a non-Hermitian matrix $A_n = \frac{1}{\sqrt{n}} M_n$ is the presence of \emph{pseudospectrum} - complex numbers $z$ for which the resolvent $(A_n - zI)^{-1} = (\frac{1}{\sqrt{n}} M_n - zI)^{-1}$ exists but is extremely large.  It is therefore of importance to obtain bounds on this resolvent, which leads one to understand for which vectors $v \in \C^n$ is $(A_n - zI) v$ likely to be small.  Expanding out the vector $(A_n - zI)v$, one encounters expressions such as $\xi_1 v_1 + \ldots + \xi_n v_n$, where $v_1,\ldots,v_n \in \C$ are fixed and $\xi_1,\ldots,\xi_n$ are iid random variables.  The problem of understanding ths distribution of such random sums is known as the \emph{Littlewood-Offord problem}, and we now pause to discuss this problem further.

\subsection{The Littlewood-Offord problem}

Let $\bv =\{v_1, \dots, v_n\}$ be a set of $n$ integers and let $\xi_1,
\dots, \xi_n$ be i.i.d random Bernoulli variables. Define $S:=
\sum_{i=1}^n \xi_i v_i$ and $p_{\bv} (a) := \P (S=a)$ and $p_{\bv}
:= \sup_{a \in \BZ} p_{\bv } (a)$.

In their study of random polynomials, Littlewood and Offord
\cite{LO} raised the question of bounding $p_{\bv}$. They showed
that if the $v_i$ are non-zero, then $p_{\bv} = O(\frac{\log
n}{\sqrt n})$. Very soon after, Erd\H {o}s \cite{Erd1}, using
Sperner's lemma, gave a beautiful combinatorial proof for the
following refinement.

\begin{theorem} \label{theorem:LO} Let $v_1, \dots, v_n$ be non-zero
numbers and $\xi_i$ be i.i.d Bernoulli random variables. Then\footnote{We use the usual asymptotic notation in this paper, thus $X = O(Y)$, $Y = \Omega(X)$, $X \ll Y$, or $Y \gg X$ denotes an estimate of the form $|X| \leq CY$ where $C$ does not depend on $n$ (but may depend on other parameters).  We also let $X = o(Y)$ denote the bound $|X| \leq c(n) Y$, where $c(n) \to 0$ as $n \to \infty$.}
$$ p_{\bv} \le \frac{\binom{n}{{\lfloor n/2 \rfloor}}}{2^n} = O(\frac{1}{\sqrt n}). $$\end{theorem}

 Notice that the bound is sharp, as can be seen from the example $\bv :=\{1,
\dots, 1\}$, in which case $S$ has a binomial distribution. Many mathematicians realized that while the classical bound in Theorem \ref{theorem:LO} is sharp as stated, it can be improved significantly under additional
assumptions on $\bv$. For instance, Erd\H {o}s and Moser \cite{EM} showed that if
the $v_i$ are distinct, then

$$p_{\bv}   =O(n^{-3/2} \ln n). $$

They conjectured that the logarithmic term is not necessary and this
was confirmed by S\'ark\"ozy and Szemer\'edi \cite{SS}. Again, the
bound is sharp (up to a constant factor), as can be seen by taking
$v_1,\ldots,v_n$ to be a proper arithmetic progression such as
$1,\ldots,n$. Stanley \cite{stanley} gave a different proof that
also classified the extremal cases.

A general picture was given by Hal\'asz, who showed, among other things,
that if one forbids more and more additive structure\footnote{Intuitively, this is because the less additive structure one has in the $v_i$, the more likely the sums $S$ are to be distinct from each other.  In the most extreme case, if the $v_i$ are linearly independent over the rationals $\Q$, then the sums $2^n$ sums $S$ are all distinct, and so $p_\bv = 1/2^n$ in this case.} in the
$v_i$, then one gets better and better bounds on $p_{\bv}$. One
corollary of his results (see \cite{Hal} or  \cite[Chapter
9]{TVbook} is the following.

\begin{theorem} \label{theorem:halasz} Consider $\bv= \{v_1, \dots, v_n\}$.
Let $R_k$ be the number of solutions to the equation

$$\eps_1 v_{i_1} + \dots  +\eps_{2k} v_{i_{2k}} =0 $$

\noindent where $\eps_i \in \{-1,1\}$ and $i_1, \dots, i_{2k} \in
\{1,2, \dots, n \}$. Then

$$p_ {\bv}= O_{k} ( n^{-2k-1/2} R_k ). $$
\end{theorem}

\begin{remark}
Several variants of Theorem \ref{theorem:LO} can be found in
\cite{Kat, GLOS, FF, Kle} and the references therein.  The connection between the Littlewood-Offord problem and random matrices was first made in \cite{KKS}, in connection with the question of determining how likely a random Bernoulli matrix was to be singular.  The paper \cite{KKS} in fact inspired much of the work of the authors described in this survey.
\end{remark}

\subsection{The inverse Littlewood-Offord problem}

Motivated by inverse theorems from additive combinatorics, in particular
Freiman's theorem (see \cite{freiman}, \cite[Chapter 5]{TVbook}) and a variant for random sums in
\cite[Theorem 5.2]{TVsing} (inspired by earlier work in \cite{KKS}), the authors \cite{TVinverse} brought a
different view to the problem. Instead of trying to improve the bound
further by imposing new assumptions, we aim to provide the full
picture by finding the underlying reason for the probability
$p_{\bv}$ to be large (e.g. larger than $n^{-A}$ for some fixed $A$).

Notice that the (multi)-set  $\bv$ has $2^n$ subsums, and $p_{\bv}
\ge n^{-C}$ mean that at least $2^n/ n^C$ among these take the same
value. This suggests that there should be very strong additive
structure in the set. In order to determine this structure, one
can study examples of $\bv$ where $p_{\bv}$ is large. For a
set $A$, we denote by $lA$ the set $lA := \{a_1+ \dots + a_l| a_i \in A
\}$. A natural example is the following.

\begin{example} Let $I=[-N,N]$ and $v_1, \dots, v_n$ be elements of
$I$. Since $S \in nI$, by the pigeon hole principle, $p_{\bv} \ge
\frac{1}{nI} = \Omega (\frac{1}{N})$.  In fact, a short
consideration yields a better bound. Notice that with probability at
least $.99$, we have $S \in 10 \sqrt n I$, thus again by the pigeonhole
principle, we have $p_{\bv}  = \Omega (\frac{1}{\sqrt n N})$. If we
set $N=n^C$ for some constant $C$, then
\begin{equation} \label{bound1} p_{\bv}  = \Omega (\frac{1}{n^{C+1/2}}).
\end{equation}
\end{example}

The next, and more general, construction comes from additive
combinatorics. A very important concept in this area is that of a \emph{generalized
arithmetic progression} (GAP). A set $Q$ of numbers is a \emph{GAP of
rank $d$} if it can be expressed as in the form
$$Q= \{a_0+ x_1a_1 + \dots +x_d a_d| M_i \le x_i \le M_i' \hbox{ for all } 1 \leq i \leq d\}$$
for some $a_0,\ldots,a_d,M_1,\ldots,M_d,M'_1,\ldots,M'_d$.

It is convenient to think of $Q$ as
the image of an integer box $B:= \{(x_1, \dots, x_d) \in \Z^d| M_i \le x_i
\le M_i' \} $ under the linear map
$$\Phi: (x_1,\dots, x_d) \mapsto a_0+ x_1a_1 + \dots + x_d a_d. $$
The numbers $a_i$ are the \emph{generators } of $P$, and $\Vol(Q) := |B|$ is the \emph{volume} of $B$. We say that $Q$ is
\emph{proper} if this map is one to one, or equivalently if $|Q| = \Vol(Q)$.  For non-proper GAPs, we of course have $|Q| < \Vol(Q)$.

\begin{example} Let $Q$ be a proper  GAP of rank $d$ and volume $V$.
Let $v_1, \dots, v_n$ be (not necessarily distinct) elements of
$P$. The random variable $S =\sum_{i=1}^n \xi_i v_i$ takes values in
the GAP $nP$. Since $|nP| \le \Vol (nB) = n^d V$, the pigeonhole
principle implies that $p_{\bv}  \ge \Omega (\frac{1}{n^d V})$. In
fact, using the same idea as in the previous example, one can
improve the bound to $\Omega (\frac{1}{n^{d/2} V})$. If we set
$N=n^C$ for some constant $C$, then
\begin{equation} \label{bound2} p_{\bv}  = \Omega (\frac{1}{n^{C+d/2}}).
\end{equation}
\end{example}

\noindent The above examples show that if the elements of $\bv$
belong to a proper GAP with small rank and small cardinality then
$p_{\bv}$ is large. A few years ago, the authors \cite{TVinverse}
showed that this is essentially the only reason:

\begin{theorem}[Weak inverse theorem]\label{theorem:weak} \cite{TVinverse} Let  $C, \epsilon > 0$ be arbitrary constants.
There are constants $d$ and $C'$ depending on $C$ and $\epsilon$
such that the following holds.
 Assume that $\bv = \{v_1, \ldots, v_n\}$ is a multiset of integers satisfying
$p_{\bv} \geq n^{-C}$. Then there is a GAP $Q$ of rank at most $d$
and volume at most $n^{C'}$ which contains all but at most
$n^{1-\epsilon}$ elements of $\bv$ (counting multiplicity).

\end{theorem}

\begin{remark}
The presence of the small set of exceptional elements is
not completely avoidable. For instance, one can add $o(\log n)$ completely
arbitrary elements to $\bv$ and only decrease $p_{\bv}$ by a factor
of $n^{-o(1)}$ at worst.  Nonetheless we expect the number of such elements to be less than what is given by the results here.
\end{remark}

The reason we call Theorem \ref{theorem:weak} {\it weak} is the fact
that the dependence between the parameters is not optimal and does
not yet reflect the relations in \eqref{bound1} and \eqref{bound2}.
Recently, we were able to modify the approach to obtain an almost
optimal result.

\begin{theorem}[Strong inverse theorem] \label{theorem:strong}
\cite{TVinversestrong}  Let $C$ and $1> \eps$ be positive constants.
Assume that
$$p_{\bv} \ge n^{-C}. $$
 Then there exists a GAP $Q$ of rank $d= O_{C, \eps} (1)$ which contains all but $O_d(n^{1 -\eps} )$
 elements of $\bv$ (counting multiplicity), where
 $$|Q| = O_{C, \eps} (n^{C - \frac{d}{2} + \eps}). $$
\end{theorem}

The bound on $|Q|$ matches \eqref{bound2}, up to the $\eps$
term. The proofs of Theorem \ref{theorem:weak} and
\ref{theorem:strong} use harmonic analysis, combined with results
from the theory of random walks and several facts about GAPs. Both
theorems hold in a more general setting, where the elements of $\bv$
are from a torsion-free group. The lower bound $n^{-C}$ on $p_{\bv}$
can also be relaxed, but the statement is more technical.

As an application of Theorem \ref{theorem:strong}, one can deduce, in
a straightforward  manner, a slightly weaker version of the {\it
forward} results mentioned above. For instance, let us show  if the
$v_i$ are different, then $p_{\bv } \le n^{-3/2+\delta}$ (for any
constant $\delta >0$). Assume otherwise and set $\eps := \delta
/2$. Theorem \ref{theorem:strong} implies that most of $\bv$ is
contained in a GAP $Q$ of rank $d$ and cardinality at most $O(
n^{3/2-\delta -d/2 + \eps }) =O(n^{1-\delta/2})=o(n)$. But since
$\bv$ has $(1-o(1))n$ elements in $Q$, we obtain a contradiction.

Next we  consider another application of Theorem
\ref{theorem:strong}, which will be more important in later
sections.  This theorem enables us execute very precise counting
arguments. Assume that we would like to count the  number of
(multi)-sets $\bv$ of integers with $\max |v_{i}| \le N$ such that
 $P(v) \ge p:= n^{-C}$.

 Fix $d \ge 1$, fix\footnote{A more detailed version of Theorems
 \ref{theorem:weak} and \ref{theorem:strong} tells us that there are
 not too many ways to choose the
 generators of $Q$. In particular, if $N =n^{O(1)}$, the number of ways to fix these is
negligible.} a  GAP $Q$  with rank $d$ and
 volume $V = n^{C -d/2}$. The dominating term will be the number of multi-subsets
of size $n$  of $Q$, which is

\begin{equation}\label{discretcounting}  |Q|^{n }= n^{(C-d/2 +\epsilon)n} \le n^{Cn} n^{-n/2+\epsilon n
}= p^{-n} n^{-n(1/2-\epsilon) }.
\end{equation}

For later purposes, we  need a continuous version of this result.
Let the $v_i$ be complex numbers. Instead of $p_{\bv}$, consider the
maximum {\it small ball} probability

$$p_{\bv}(\beta)  =\max_{z \in \C} \P (|S-z| \le \beta) . $$

 Given a small $\beta >0$ and $ p= n^{-O(1)}$, the
collection of $\bv$ such that $|v|=1$ and $p_{\bv}(\beta) \ge p$ is
infinite, but we are able to show that it can be approximated by a
small set.

\begin{theorem} [The $\beta$-net Theorem] \cite{TVcir1}  Suppose that $p = n^{-O(1)}$. Then the set
of unit vectors $\bv= (v_{1}, \dots, v_{n})$ such that
$p_{\bv}(\beta)  \ge p$ admits an  $\beta$-net (in the infinity
norm\footnote{In other words, for any $\bv$ with $p_{\bv}(\beta) \geq p$, there exists $\bv' \in \Omega$ such that all coefficients of $\bv - \bv'$ do not exceed $\beta$ in magnitude.} $\Omega$ of size at most

\begin{equation} \label{contcounting} |\Omega| \leq p^{-n} n^{-n/2 +o(n)} .
\end{equation} \end{theorem}

\begin{remark}\label{Rvrem} A related result (with different parameters) appears in \cite{RV}; in our notation, the probability $p$ is allowed to be much smaller, but the net is coarser (essentially, a $\beta \sqrt{n}$-net rather than a $\beta$-net).  In terms of random matrices, the results in \cite{RV} are better suited to control the extreme tail of such quantities as the least singular value of $A_n - zI$, but require more boundedness conditions on the matrix $A_n$ (and in particular, bounded operator norm) due to the coarser nature of the net.
\end{remark}

\section{Computer Science}

Our next stop is computer science and numerical linear algebra, and in particular the problem of dealing with \emph{ill-conditioned} matrices, which is closely related to the issue of pseudospectrum which is of central importance in the circular law problem.

\subsection{Theory vs Practice}

Running times of algorithms are frequently estimated by worst-case
analysis. But in practice, it has been observed that many
algorithms, especially those involving a large matrix,  perform
significantly better than the worst-case scenario. The most
famous example is perhaps the simplex algorithm in linear
programming. Here, the basic problem (in its simplest form) is to
optimize a goal function $c \cdot x$, under the constraint $Ax \le
b$, where $c, b$ are given vectors of length $n$ and $A$ is an $n
\times n $ matrix. In the worst case scenario, this algorithm
takes exponential time. But in practice, the algorithm runs
extremally well. It is still very popular today, despite the fact
that there are many other algorithms proven to have polynomial
complexity.

There have been various attempts to explain this phenomenon. In this section we will discuss an influential recent explanation given by Spielman and Teng \cite{ST,
ST1}.

\subsection {The effect of noise}
 An important issue in the theory of computing is noise,
as almost all  computational processes are
 effected by it.  By the word \emph{noise}, we would like to  refer to all kinds of
 errors occurring in a process, due to both humans and machines, including
 errors  in measuring,  errors caused by truncations,
errors committed in transmitting  and inputting the data, etc.

Spielman and Teng \cite{ST} pointed out that  when we are
interested in a solving a certain system of
 equations, because of noise,
 our computer actually ends
 up solving a slightly perturbed version of the system. This is the
 core of their so-called {\it smooth analysis} that they used to
 explain the effectiveness of a specific algorithm (such as the simplex
 method). Interestingly, noise, usually a burden, plays a ``positive''
 role here, as it smoothes the inputs randomly, and so prevents a very
 bad input from ever occurring.

  The puzzling question here is, of
course: {\it why is the perturbed input typically better than the
original (worst-case) input ?}

In order to give a mathematical explanation, we  need to introduce
some notion. For an $n \times n$ matrix $M$, the \emph{condition
number} $\kappa(M)$ is defined as
$$\kappa(M):= \|M\| \| M^{-1} \|$$
where $\| \|$ denotes the operator norm.
(If $M$ is not invertible, we set $\kappa (M) =\infty$.)

The condition number  plays a crucial role in numerical linear algebra; in particular, the condition number $\kappa(M)$ of a matrix $M$  serves as a simplified proxy for the accuracy and stability of most algorithms used to solve the
 equation $Mx=b$ (see \cite{BT, GvL}, for example). The
 exact solution $x= M^{-1} b$, in theory, can be computed quickly (by
 Gaussian elimination, say). However, in practice computers can only represent a
 finite subset of real numbers and this leads to two
 difficulties: the represented numbers cannot be arbitrarily large
 or small, and there are gaps between two adjacent represented numbers. A quantity which is frequently used in numerical
 analysis is  $\eps_{\mathrm{machine}}$ which is half of the distance
 from $1$ to the nearest represented number. A fundamental result
 in numerical analysis \cite{BT}
 asserts that if one denotes by $\tilde x$ the result computed by
 computers, then the relative error $\frac{ \| \tilde x - x \|
 }{\|x\|}$ satisfies

 $$ \frac{ \| \tilde x - x \| }{\|x\|} = O\big( \eps_{\mathrm{machine}}
 \kappa(M) \big) $$

Following the literature, we call $M$ {\it well-conditioned} if
$\kappa (M)$ is small. For quantitative purposes, we say that an $n$
by $n$ matrix $M$ is {\it well-conditioned} if its condition number
is polynomially bounded in $n$ (that is, $\kappa(M) \le n^C$ for some
constant $C$ independent of $n$).

\subsection{Randomly perturbed matrices are well-conditioned}

The analysis in \cite{ST} is guided by the following fundamental intuition\footnote{This conjecture, of course, does not fully explain the phenomenon of smoothed analysis, since it may be that a well-conditioned matrix still causes a difficulty in one's linear algorithms for some other reason, or perhaps the original ill-conditioned matrix did not cause a difficulty in the first place; we thank Alan Edelman for pointing out this subtlety.  Nevertheless, Conjecture \ref{con} does provide an informal intuitive justification of smoothed analysis, and various rigorous versions of this conjecture were used in the formal arguments in \cite{ST}: see Section 1.4 of that paper for further discussion.}:

 \begin{conjecture}\label{con}  For every input instance, it is unlikely
that a slight random perturbation of that instance has large
condition number. \end{conjecture}

More quantitatively,

\begin{conjecture} Let $A$ be an arbitrary $n$ by $n$ matrix and let $M_n$ be a random
 $n$ by $n$ matrix.
 Then with high probability $A+M_n$ is well-conditioned.
 \end{conjecture}

\vskip2mm Notice that here one allows  $A$ to have a  large
condition number.

Let us take a look at $\kappa (A+M_n) = \| A+M_n \| \| (A+M_n)^{-1}
\|$. In order to have $\kappa (A+M_n) = n^{O(1)}$, we want to upper-bound
both $\| A+M_n \| $ and $\| (A+M_n)^{-1} \|$. Bounding $\| A+M_n \|$ is
easy, since by the triangle inequality
$$ \|A+M_n \| \le \|A\| + \|M_n \|. $$

In most models of random matrices, $\|M_n \| \leq n^{O(1)}$ with very
high probability, so it suffices to assume that $\|A \| \leq n^{O(1)}$; thus we assume that the matrix $A$ is of polynomial size compared to the noise level.  This is a fairly reasonable assumption for high-dimensional matrices for which the effect of noise is non-negligible\footnote{In particular, it is naturally associated to the concept of \emph{polynomially smoothed analysis} from \cite{ST}.}, and we
are going to assume it in the rest of this section. 

The remaining problem is to bound the norm of the inverse $\|
(A+M_n)^{-1} \|$. An important detail here is how to choose the
random matrix $M_n$. In their works \cite{ST, ST1, SST}, Spielman
and Teng (and coauthors) set $M_n$ to have iid Gaussian entries
(with variance 1) and obtained the following bound, which played a
critical role in their smooth analysis \cite{ST, ST1}.

\begin{theorem} \label{theorem:STcondition} Let $A$ be an arbitrary $n$ by $n$ matrix
and $M_n$ be a random matrix with iid Gaussian entries.
 Then for any $ x >0$,
$$\P( \|(A+M_n)^{-1} \|  \ge x) = O(\frac{\sqrt n}{x} ) . $$
\end{theorem}

While Spielman-Teng smooth analysis does seem to have the right
philosophy, the choice of $M_n$ is a bit artificial. Of course, the
analysis still passes if one replaces Gaussian by a fine enough
approximation.   A large fraction of problems in linear programming
deal with integral matrices, so the noise is perturbation by integers. In other
cases, even when the noise has continuous support, the data is
strongly truncated.  For example, in many engineering problems, one
does not keep more than, say,  three to   five decimal places. Thus, in many
situations, the entries of $M_n$ end up having discrete support
with relatively small size, which may not even grow with $n$, while
the approximation mentioned above would require this support to have
size exponential in $n$. Therefore, in order to come up with an
analysis that better captures real life data, one needs to come up
with a variant of Theorem \ref{theorem:STcondition} where the
entries of $M_n$ have discrete support.

This problem was suggested to the authors by Spielman a few years ago. Using the
Weak Inverse Theorem, we were able to prove the following variant
of Theorem \ref{theorem:STcondition} \cite{TVstoc}.

\begin{theorem} \label{theorem:conditionTV} For any constants $a,c >0$, there is a constant
$b=b(a,c)>0$ such that the following holds.
 Let $A$ be an $n$ by $n$
matrix such that $\|A\|\le  n^{a}$  and let $M_n$ be a random matrix
with iid Bernoulli entries.
 Then
$$\P( \|(A+M_n)^{-1} \|  \ge n^b)\le n^{-c }. $$
\end{theorem}

Using the stronger $\beta$-net  Theorem, one can have a nearly optimal
relation between the constants $a$, $b$ and $c$ \cite{TVgeneral}. These
results extend, with the same proof, to a large variety of
distributions. For example, one does not need require the entries of
$M_n$ to be iid\footnote{In practice, one would expect the noise at a large
entry to have larger variance than one at a small entry, due to multiplicative effects.}, although independence is crucially exploited in the proofs. Also,
one can allow many of the entries to be 0 \cite{TVstoc}.

\begin{remark}  Results of this type first appear in \cite{Rud} (see also \cite{Lit} for some earlier related work for the least singualar value of \emph{rectangular} matrices).  In the special case where $A=0$ and where the entries of $M_n$ are iid and have finite fourth moment, Rudelson and Vershynin \cite{RV}  (see also \cite{RV2}, \cite{RV3}) obtained sharp
 bounds for $\|(A+M_n)^{-1}\|$, using a somewhat different method, which relies on an inverse
 theorem of a slightly different nature; see Remark \ref{Rvrem}.
\end{remark}

The main idea behind the proof of Theorem \ref{theorem:conditionTV}, which first appears in \cite{Rud}, is the following. Let $d_{i}$ be the distance from the $i^{\operatorname{th}}$ row
vector of $A+M_n$  to the subspace spanned by the rest of the rows.  Elementary linear algebra (see also \eqref{neg} below) then gives the bound
$$\| (A+M_n)^{-1} \| =n^{O(1)}  (\min_{1 \leq i \leq n} d_{i} )^{-1}.
$$
Ignoring various factors of $n^{O(1)}$, the main task is then to understand the distribution of $d_i$ for any given $i$.

If  $v= (v_{1}, \dots, v_{n})$ is the normal vector of a hyperplane $V$,
then the distance from a random vector $(a_1+ \xi_1, \dots, a_n+
\xi_n)$ to the hyperplane $V$ is given by the formula
$$ | v_{1 } (\xi_{1 }+a_1)  + \dots + v_{n} (\xi _{n }+a_n)  | =|\sum_i a_i v_i + S |
$$
where $S := \sum_{i=1}^n v_i \xi_i$ is as in the previous section.

To estimate the chance that $|\sum_{i=1}^n a_i v_i + S| \le \beta$,
 the notion of the small ball probability $p_{\bv}(\beta)$ comes naturally. Of course, this quantity depends on the normal vector $\bv$, and so we now divide into cases depending on the nature of this vector.

If $p_{\bv}(\beta)$ small, we can be done using a conditioning
argument\footnote{Intuitively, the idea of this conditioning argument is to first fix (or ``condition'') on $n-1$ of the rows of $A+M_n$, which should then fix the normal vector $\bv$.  The remaining row is independent of the other $n-1$ rows, and so should have a probability at most $p_\bv(\beta)$ of lying within $\beta$ of the span of the those rows.  There are some minor technical issues in making this argument (which essentially dates back to \cite{Komlos}) rigorous, arising from the fact that the $n-1$ rows may be too degenerate to accurately control $\bv$, but these difficulties can be dealt with, especially if one is willing to lose factors of $n^{O(1)}$ in various places.}. On the other hand, the $\beta$-net Theorem says
that there are ``few'' $\bv$ such that  $p_{\bv}(\beta)$ is large, and in this case a direct counting argument
finishes the job\footnote{For instance, one important class of $\bv$ for which $p_\bv(\beta)$ tends to be large are the \emph{compressible} vectors $\bv$, in which most of the entries are close to zero.  Each compressible $\bv$ (e.g. $\bv = (1,-1,0,\ldots,0)$) has a moderately large probability of being close to a normal vector for $A+M_n$ (e.g. in the random Bernoulli case, $\bv = (1,-1,0,\ldots,0)$ has a probability about $2^{-n}$ of being a normal vector); but the number (or more precisely, the metric entropy) of the set of compressible vectors is small (of size $2^{o(n)}$) and so the net contribution of these vectors is then manageable.  Similar arguments (relying heavily on the $\beta$-net theorem) handle other cases when $\bv$ is large (e.g. if most entries of $\bv$ live near a GAP of controlled size).}.  Details can be found in \cite{TVstoc}, \cite{TVcir1}, or \cite{TVgeneral}.

\section{Back to probability}

\subsection{The replacement principle}

Let us now take another look at the Circular Law Conjecture. Recall
that $\lambda_{1}, \dots, \lambda_{n}$ are the eigenvalues of
$A_n = \frac{1}{\sqrt n} M_{n}$, which generates a normalized counting measure
$\mu_{A_n}$. We want to show that $\mu_{A_n}$ tends (in probability) to the
uniform measure $\mu$ on the unit disk.

The traditional way to attack this conjecture is via a Stieltjes
transform technique\footnote{The more classical \emph{moment method}, which is highly successful in the Hermitian setting (for instance in proving Theorem \ref{theorem:Pastur}), is not particularly effective in the non-Hermitian setting, because moments such as $\tr A_n^m$ for $m=0,1,2,\ldots$ do not determine the ESD $\mu_{A_n}$ (even approximately) unless one takes $m$ to be as large as $n$; see \cite{bai}, \cite{BS} for further discussion.}, following \cite{Girko1, bai}.  Given a
(complex) measure $\nu$, define, for any $z$ with Im $z >0$,
$$s_{\nu}(z) :=  \int \frac{1}{x-z} d \nu (x). $$
For the ESD $\mu_{A_n}$, we have
$$s_{\mu_{A_n}} (z) = \frac{1}{n} \sum \frac{1}{\lambda_{i} - z } .
$$

Thanks to standard results from probability\footnote{One can also use the theory of logarithmic potentials for this, as is done for instance in \cite{GT1}, \cite{PZ}.}, in order to establish
the Circular Law Conjecture in the strong (resp. weak) sense, it suffices to show that $s_{\mu_n}(z)$
converges almost surely (resp. in probability) to $s_{\mu}(z)$ for almost all $z$ (see \cite{TVcir2} for a precise statement).

Set $z=: s+ it$ and $s_{n} (z) =: S+ iT$. Since $s_n$ is analytic
except at the poles, and vanishes at infinity, the Stieltjes transform $s_n(z)$ is determined by its the real part $S$. Let us
take a closer look at this variable:

\begin{eqnarray*} S &=& \frac{1}{n} \sum \frac{ \Re(\lambda_{i}) -s} {| \lambda_{i} -z |^{2} } \\
&=&- \frac{1}{2n} \sum _{} \frac{\partial}{ \partial s}  \log | \lambda_{i} -z |^{2}  \\
&=& - \frac{1}{2} \frac{\partial}{\partial s} \int_{0}^{\infty}
\log x  \,\,  \partial \eta_{n}
\end{eqnarray*}

where
$$\eta_{n} := \mu_{(\frac{1}{\sqrt n} M_n - zI)(\frac{1}{\sqrt n} M_n - zI)^*}$$
is the normalised counting measure of the (squares of the)
\emph{singular values} of $\frac{1}{\sqrt n} M_{n} -zI$. Notice that
in the third equality, we use the fact that $\prod |\lambda_{i}-z| =
|\det (\frac{1}{\sqrt n} M_{n} - zI) |$. This step is critical as it
reduces the study of a complex measure to a real one, or in other words to study the ESD of a Hermitian matrix rather than a non-Hermitian matrix.

Putting this observation in the more general setting of Theorem
\ref{theorem:main1}, we arrived at the following useful result.

\begin{theorem}[Replacement principle]\label{theorem:replacement}\cite{TVcir2} Suppose for each $n$
that $A_n, B_n \in M_n(\BBC)$ are ensembles of random matrices.
  Assume that
\begin{itemize}
\item[(i)] The expression
\begin{equation}\label{pan}
\frac{1}{n^2} \|A_n\|_F^2 + \frac{1}{n^2} \|B_n\|_F^2
\end{equation}
is weakly (resp. strongly) bounded\footnote{A sequence $x_n$ of non-negative random variables is said to be \emph{weakly bounded} if $\lim_{C \to \infty} \liminf_{n \to \infty} \P( x_n \leq C ) = 1$, and \emph{strongly bounded} if $\limsup_{n \to \infty} x_n < \infty$ with probability $1$.}
\item[(ii)] For almost all complex numbers
$z$, $$\frac{1}{n}
 \log |\det(\frac{1}{\sqrt{n}} A_n - zI)| - \frac{1}{n} \log |\det(\frac{1}{\sqrt{n}} B_n -
 zI)|$$
converges weakly (resp. strongly) to zero.  In particular, for each
fixed $z$, these determinants are non-zero with probability $1-o(1)$
for all $n$ (resp. almost surely non-zero for all but finitely many
$n$).
\end{itemize}
Then $\mu_{\frac{1}{\sqrt{n}} A_n} - \mu_{\frac{1}{\sqrt{n}} B_n}$
converges weakly (resp. strongly) to zero.
\end{theorem}

At a technical level, this theorem reduces Theorem \ref{theorem:main1} to the
comparison of $\log |\det(\frac{1}{\sqrt{n}} A_n - zI)| $ and $\log
|\det(\frac{1}{\sqrt{n}} B_n -
 zI)|$.

\begin{remark}
Note that this expression is large and unstable when $z$ lies in the \emph{pseudospectra} of either $\frac{1}{\sqrt{n}} A_n$ or $\frac{1}{\sqrt{n}} B_n$, which means that the resolvent $(\frac{1}{\sqrt{n}} A_n - zI)^{-1}$ or $(\frac{1}{\sqrt{n}} B_n - zI)^{-1}$ is large.  Controlling the probability of the event that $z$ lies in the pseudospectrum is therefore an important portion of the analysis.  This technical problem is not an artefact of the method, but is in fact essential to any attempt to control non-Hermitian ESDs for general random matrix models, as such ESDs are extremely sensitive to perturbations in the matrix in regions of pseudospectrum.  See \cite{bai}, \cite{BS} for further discussion.
\end{remark}

\subsection{Treatment of the pole}

Using techniques from probability, such as the moment method, one
can show that the distributions of the singular values of
$\frac{1}{\sqrt{n}} A_n - zI$ and $\frac{1}{\sqrt{n}} B_n - zI$ are
asymptotically the same\footnote{In the setting where the matrices $X_n$ and $Y_n$ have iid entries, one can use the results of \cite{doz} to establish this.  In the non-iid case, an invariance principle from \cite{chat} gives a slightly weaker version of this equivalence; this was observed by Manjunath Krishnapur and appears as an appendix to \cite{TVcir2}.} \cite{bai, TVcir1, doz, TVcir2, chat}. This,
however, is not sufficient to conclude that $\frac{1}{n} \log
|\det(\frac{1}{\sqrt{n}} A_n - zI)| $ and $\frac{1}{n} \log
|\det(\frac{1}{\sqrt{n}} B_n -
 zI)|$ are close.  As remarked earlier, the main difficulty here is that some of the singular values can be
very small and thus significantly influence the value of logarithm.

Now is where Theorem \ref{theorem:conditionTV} enters the
picture. This theorem  tells us that (with overwhelming
probability), there is no mass between $0$ and (say) $n^{-C}$, for
some sufficiently large constant $C$.  Using this critical
information, with some more work\footnote{In particular, the presence of certain factors of $\log n$ arising from inserting Theorem \ref{theorem:conditionTV} into the normalized log-determinant $\frac{1}{n} \log |\det(\frac{1}{\sqrt{n}} A_n - zI)|$ forces one to establish a \emph{convergence rate} for the ESD of $\frac{1}{\sqrt{n}} A_n - zI$ which is faster than logarithmic in $n$ in a certain sense.  This is what ultimately forces one to assume the bounded $(2+\eta)^{\operatorname{th}}$ moment hypothesis.  Actually the method allows one to relax this hypothesis to that of assuming $\E |\a|^2 \log^C (2+|\a|) < \infty$ for some absolute constant $C$ (e.g. $C=16$ will do).}, we obtain:

\begin{theorem} \label{theorem:weakCL} \cite{TVcir1}   The Circular Law holds (with both strong and weak convergence) under
the extra condition that the entries have bounded $(2+\eta)^{\operatorname{th}}$ moment, for
some constant $\eta >0$. \end{theorem}

\begin{remark}
Shortly after the appearance of \cite{TVcir1}, G\"otze and Tikhomirov \cite{GT2} gave an alternate proof of the weak circular law with these hypothesis, using a variant of Theorem \ref{theorem:conditionTV},
which they obtained via a method from \cite{Rud}, \cite{RV}. This method is
based on a different version of the Weak Inverse Theorem.
\end{remark}

\subsection{Negative second moment and sharp concentration}

At the point it was written, the analysis in \cite{TVcir1} looked
close to the limit of the method. It took some time to realize where the extra
moment condition came from and even more time to figure out a way to avoid
that extra condition. Consider the sums

$$\frac{1}{n} \log |\det(\frac{1}{\sqrt{n}} A_n - zI)|  = \frac{1}{n} \sum_{i=1}^n \log \sigma_i, $$
where $\sigma_1 \ge \dots \ge \sigma_n$ are the
singular values of $\frac{1}{\sqrt n} A_n -zI$, and
$$\frac{1}{n} \log |\det(\frac{1}{\sqrt{n}} B_n - zI)| = \frac{1}{n} \sum_{i=1}^n \log \sigma'_i, $$
where $\sigma'_1 \ge \dots \ge \sigma'_n$ are the
singular values of $\frac{1}{\sqrt n} B_n -zI$.

As already mentioned, we know that the bulk of the $\sigma_i$ and $\sigma_i'$
are distributed similarly. For the smallest few, we used the lower
bound on $\sigma_n$ as a uniform bound be show that their
contribution is negligible. This turned out to be wasteful, and we
needed to use the extra moment assumption to compensate the loss in
this step.

In order to remove this assumption, we need to find a way to give
a better bound on other singular values. An important first step is
the discovery of the following simple, but useful, identity.

{\bf The Negative Second Moment Identity.} \cite{TVcir2}  Let $A$ be an
$m \times n$ matrix, $m \le n$. Then
\begin{equation}\label{neg}
\sum_{i=1} ^{m}  d_{i}^{-2} = \sum_{i=1} ^{m} \sigma_{i} ^{-2}
\end{equation}
where, as usual,  $d_{i} $ are the distances and $\sigma_{i} $ are
the singular values.

One can prove this identity  using undergraduate linear algebra.
With this in hand, the rest of the proof falls into place\footnote{A possible alternate approach would be to bound the intermediate singular values directly, by adapting the results from \cite{RV2}.  This would however require some additional effort; for instance, the results in \cite{RV2} assume zero mean and bounded operator norm, which is not true in general when considering $\frac{1}{\sqrt{n}} A_n - zI$ for non-zero $z$ assuming only a mean and variance condition on the entries of $A_n$.  In any case, the analysis in \cite{RV2} ultimately goes through a computation of the distances $d_i$, similarly to the approach we present here based on the negative second moment identity.}.  Consider
the singular values $\sigma_{1} \ge \dots \ge \sigma_{n} $ involved
in our analysis, and use $A$ as shorthand for $\frac{1}{\sqrt n} A_n
-zI$.
 To bound $\sigma_{n-k}$ from below, notice that  by the
interlacing law
$$\sigma_{n-k} (A) \ge \sigma_{m-k} (A' )$$

where $m:=n-k$ and $A'$ is an $m \times n $ truncation of $A$,
obtained by omitting the last $k$ rows. The  Negative Second Moment
Identity implies

$$k \sigma_{m-k} (A') ^{{-2} } \le \sum _{i=1}^{m} \sigma_{i} (A')^{-2} = \sum_{i=1}^{m} d_{i} ^{-2}  . $$

On the other hand, the right-hand side can be bounded efficiently,
thanks to the fact that all $d_i$ are large with overwhelming
probability, which, in turn, is a consequence of Talagrand's
inequality \cite{Tal}:

\begin{lemma}[Distance Lemma]\cite{TVdet, TVcir2} With probability
 $1- n^{-\omega(1)}$, the distance from a random
row vector to a subspace of co-dimension $k$ is  at least
$\frac{1}{100} \sqrt {k/n}$, as long as $k \gg {\log n }$.
\end{lemma}

Thus, with overwhelming probability, $ \sum_{i=1}^{m}d_{i}^{-2}$ is
$ \Omega (m/nk)= \Omega((n-k)/nk)$, which implies

{$$\sigma_{n-k}(A) \ge \sigma_{m-k} (A') \gg \frac{k}{\sqrt{(n-k)n}}. $$}

This lower bound now is sufficient to establish Theorem
\ref{theorem:main1} and with it the Circular Law in full generality.

\section{Open problems}

Our investigation leads to open problems in several areas:

\vskip2mm

{\it Combinatorics.}  Our studies of Littewood-Offord problem focus on the  linear form
$S:=\sum_{i=1}^{n } v_{i }xi_{i} $. What can one say about higher degree polynomials ?

In \cite{CTV}, it was shown that for a quadratic form $Q:=\sum_{1\le i,j \le n} c_{ij}\xi_{i}\xi_{j}$ with non-zero coefficients, $\P(Q=z)$ is  $O(n^{-1/8})$. It is simple to improve this bound to
$O(n^{-1/4} )$ \cite{CV1}. On the other hand, we conjecture that the truth is
$O(n^{-1/2} )$, which would be sharp by taking $Q= (\xi_{1} + \dots + \xi_{n} )^{2} $.
Costello (personal communication) recently improved the bound to $O(n^{-3/8})$, and it looks likely that his approach will lead to the optimal bound, or something close.

The situation with higher degrees is much less clear. In \cite{CTV}, a bound of the form
$O(n^{-c_{k} } )$ was shown, where $c_{k} $ is a positive constant depending on $k$, the degree of the polynomial involved. In this bound $c_{k}$ decreases very fast with $k$.

\vskip2mm

{\it Smooth analysis.}  Spielman-Teng smooth analysis of the simplex algorithm \cite{ST} was done with gaussian noise. It is a very interesting problem to see if one can achieve the same conclusion with
discrete noise with fixed support, such as Bernoulli. It would give an even more convincing explanation to the efficiency
of the simplex method. As discussed earlier, noise that occurs in practice typically has discrete, small support. (This question was mentioned to us by several researchers, including Spielman,  few years ago.)

As discussed earlier, we now have the discrete version of Theorem \ref{theorem:STcondition}. While
Theorem \ref{theorem:STcondition} plays a very important part in Spielman-Teng analysis \cite{ST1},
there are several other parts of the proof that make use of the continuity of the support in
subtle ways.  It is possible to modify these parts to work for  fine enough discrete approximations  of
the continuous (noise) variables in question. However, to do so it seems one need to
make the size of the support very large (typically exponential in $n$, the size of the matrix).

Another exciting direction is to consider
even more realistic models of noise. For instance,

\begin{itemize}

\item In several problems, the  matrix may have  many {\it frozen } entries, namely those which
are not effected by noise. In particular,  an entry which is zero (by nature of the problem) is likely to
stay zero in the whole computation. It is clear that the  {\it pattern} of the frozen entries will
be of importance. For example, if the first column  consists of (frozen) zero, then no matter how
the noise effects the rest of the matrix, it will always be non-singular (and of  course
ill-conditioned). We hope to classify all
patterns where theorems such as Theorem \ref{theorem:main1} are still valid.

\item In non-frozen places, the noise could have different distributions. It is natural to think that the
error at  a large entry should have larger variance than the one occurring at a smaller entry.

\end{itemize}

Some preliminary results in these directions are obtained in \cite{TVstoc}. However, we are still at the very beginning
of the road and much needs to be done.

\vskip2mm

{\it Circular Law.} A natural question here is to investigate the rate of convergence. In \cite{TVcir1},
we observed that under the extra assumption that the $(2+\eps)$-moment of the entries are bounded,
we can have rate of convergence of order $n^{-\delta} $, for some
positive constant  $\delta $ depending on $\eps$.  The exact dependence between $\eps$ and $\delta$
is not clear.

Another question concerns the determinant  of random matrices. It is known, and not hard to prove, that
$\log |\det M_{n} |$ satisfies a central limit theorem, when the entries of $M_{n} $ are iid gaussian, see \cite{Girkodet1, CV2}. Girko \cite{Girkodet1} claimed that the same result holds for
much more general models of matrices. We, however, are unable to verify his arguments.
It would be nice to have an alternative proof.

\section{Acknowledgements}

Thanks to Peter Forrester, Kenny Maples, and especially Alan Edelman for corrections.

\end{document}